
\input amstex.tex
\documentstyle{amsppt}
\input amsppt1.tex
\nologo
\magnification=1200
\nopagenumbers
\parindent=0pt
\NoBlackBoxes
\define\op#1{\operatorname{#1}}

\topmatter
\author   A. Kazarnovski-Krol        \endauthor
\title Cycles for asymptotic solutions and Weyl group \endtitle
\address{
  Department of Mathematics
  Rutgers University
  New Brunswick, NJ 08854, USA}

\vskip 1cm

\abstract{In this paper cycles for asymptotic solutions are
described. Cycles are enumerated by elements of symmetric group.
Leading asymptotic and leading coefficient are calculated.
Value of certain multiple integral over special cycles is calculated
with the help of result of Opdam. }
\endtopmatter

\document

\head 0. Introduction \endhead
In this paper cycles for asymptotic solutions are
described. Cycles are enumerated by elements of symmetric group.
Leading asymptotic and leading coefficient are calculated.
Value of certain multiple integral over special cycles is calculated
with the help of result of Opdam.

In [33] homological treatment of Harish-Chandra
decomposition is given.

We start with recalling the necessary material from [18, 19 ].
We restrict ourselves to the case of root system of type $A_n$.
For better exposition the reader should consult directly with the papers
[18,19,38, 39, 40].
\subhead{0.1. Differential operator of second order}\endsubhead
Let $L$ be the following differential operator
$$
L=
{\sum_{i=1}^{n+1} (z_i \frac{\partial} {\partial z_i})^2} -
{k\sum_{i < j} {\frac{z_j +z_i}{z_j-z_i}}(z_i \frac{\partial} {\partial z_i}-
z_j \frac{\partial} {\partial z_j})}.
$$

Let $u_i= \log( z_i)$ for $i=1, 2, \ldots, n+1$.
In coordinates $u_1, \ldots, u_{n+1}$
operator $L$ reads as:

$$
L=
\sum_{i=1}^{n+1} \frac {\partial^2}{{\partial u_i^2}}-
k\sum_{i<j} \frac{1+e^{u_i-u_j}}{1-e^{u_i -u_j}}(\frac {\partial} {\partial u_i
}-
\frac {\partial}{\partial u_j})
$$
\remark{Remark }
In notations of Harish-Chandra cf.[19] operator $L$ reads as:
$$L=(H_1^2 +\ldots +H_l^2) +\sum_{\alpha \in P_{+}}
{{(e^{\alpha}+e^{-\alpha})}\over
{(e^{\alpha}-e^{-\alpha})}}  H_{\alpha}.$$

Operator $L$ plays the crucial role in the theory of zonal spherical
functions on noncompact Rimannian spaces $G / K$.

We, following [18], use $\alpha$ twice as usual restricted roots  and $k$
half the usual multiplicities.

The theory of zonal spherical functions goes back to H.Weyl and
E.Cartan [ 42 ]. They used mostly integral methods.  In [41 ] I.Gelfand
suggested to use  Laplace-Casimir operators , cf. also [26 ] .
\smallskip
$L$ is the radial part of Laplace-Casimir operator of second order
with respect to Cartan decomposition ($G=KAK$).

\endremark

\subhead {0.2. Root system of type $A_n$}\endsubhead

Let $V$ be a Euclidean $(n+1)$-dimensional space with inner product
$(.,.)$.

For $\alpha\in V$

$$
r_{\alpha}(\lambda)= {
\lambda - { (\lambda, \alpha^\vee)\alpha}}
$$

$(\alpha^\vee ={ 2 \alpha \over (\alpha, \alpha)} )$

denotes the orthogonal reflection with respect to the hyperplane
 $(\alpha,.)=0$.

 Let $$e_1, e_2, \ldots, e_{n+1}$$ be the orthonormal
basis of $V$.
Let $E$ be  a $n$-dimensional subspace of $V$ , which is
orthogonal to $e_1+e_2 +\ldots + e_{n+1}$ :
$$E =\{(u_1,u_2, \dots,u_{n+1} );\quad  u_1+u_2 +\ldots+u_{n+1}=0\}.$$
Then $E$ is also Euclidean and we will sometimes identify it with its dual
without additional comments.
Then the set $R$ of $n(n+1)$ vectors in $E$:
$$ R:=\{e_{i} -e_{j} | \quad i\ne{j}\}$$

is  called {\bf the root system of type $A_n$}.

\smallskip
Take a vector $\alpha \in R$, say $\alpha =e_i-e_j$.
Then reflection $r_{\alpha}$
permutes $i$th and $j$th coordinate:
$$
\align
&{r_{\alpha}}{(\lambda_1,\; \lambda_2,\ldots,\lambda_i,
\ldots,\lambda_j,\ldots,\lambda_{n+1})}\\
&\qquad\qquad = {
 (\lambda_1,\; \lambda_2,\ldots,\lambda_j,\ldots,
 \lambda_i,\ldots,\lambda_{n+1})
 }
\endalign
$$
 The Weyl group $W$ is the subgroup of isometries of $E$, which is generated by
reflections
$r_{\alpha}, \quad\alpha \in R$. Obviously, for the root system of
type $A_n$  $ W $  coincides   with
the group of permutation of $n+1$ letters $S_{n+1}$.

The set of vectors $R_{+}= \{ e_i -e_j| \quad i<j \} $ is the usual
choice for the positive roots.
The halfsum of positive roots is equal to:
$$
\delta =\frac{1}{2} \sum_{\alpha \in R_{+}}{ \alpha} =  (\frac{n}{ 2} ,\frac{
n-2}{ 2}, \ldots,
\frac{-n}{ 2})
$$
 Set
$$
\rho =\rho(k) =\frac{1}{2} k\sum_{\alpha \in R_{+}}{ \alpha} = k \delta=
 (\frac{n}{ 2}k ,\frac{ n-2}{ 2}k, \ldots,
\frac{-n}{ 2}k)   \quad .
$$

Let $\alpha_i= e_i-e_{i+1}$, $i=1,2, \ldots, n$. Then vectors
$\alpha_1, \alpha_2, \ldots, \alpha_n$ are called simple roots, and
any positive root is expressed as a linear combination of simple roots
with nonnegative coefficients.
Introduce a partial ordering on $E$ by
$$
\lambda \ge \mu\; \text{if and only if}\;
\lambda -\mu=\sum_{j=1}^{n} {l_{j} \alpha_{j}} \;\text{with}\;
{l_j \in \Bbb Z_{+}},\quad \text{where } \;\Bbb Z_{+}= 0,1,2,\ldots
$$

Let $\Lambda_1,\Lambda_2, \ldots, \Lambda_n$ be fundamental weights,
i.e.
 $\Lambda_i \in E, (\Lambda_i, \alpha_j)= \delta_{ij}$
, where $ \delta_{ij}$ is Kronecker's delta.


\subhead{0.3. Ring of differential operators  } \endsubhead
Let $\Cal{R}$ be the algebra of functions  generated by the functions
$$
\{
\frac
{1}
{1 - e^{u_i-u_j} }\;, i<j \}
$$

Consider $\Cal D$ the set of all
differential operators
with coefficients in ${\Cal{R}}$, i.e. of the form

$$
\sum_{i} f_{i}(u_1,\ldots, u_{n+1}) P_{i}\big({\frac{\partial}
{\partial{u_1}}},
\ldots,{\frac{\partial} {\partial {u_{n+1}}}}\big)
\;,
$$
where  $f_i(u_1,\ldots,u_{n+1}) \in \Cal R$.
Using the identity
$$
\frac
{1}
{1 -{ e^{u_i-u_j}}}
=
\sum_{m=0}^{\infty} {e^{m(u_i-u_j)}}
$$
For each differential operator from $\Cal D$ one obtains  a representation of
the
form:

$$
P= \sum_{\mu \ge 0} {e^{\mu(u)}
} p_{\mu}\big(\frac{\partial} {\partial u_1},\ldots, \frac {\partial}
{\partial u_{n+1}}\big) \; ,
$$
where $\mu(u)= (\mu, u)$ and $u=(u_1,\ldots, u_{n+1})$

In particular, $L$ is a differential operator of $\Cal D$  and $L$
has the following asymptotic  expansion:

$$
L=
\sum {\frac {\partial^2}{\partial u_i^2} } -2\sum_i (\rho,
e_i) \frac{\partial}{\partial u_i}-
2\sum_{i < j } k \sum_{m=1}^{\infty}  e^{m(u_i-u_j)}
(\frac {\partial}{\partial u_i}-\frac {\partial}{\partial u_j}).
$$

\smallskip
\definition{0.4. Definition }(Harish-Chandra homomorphism)
Harish-Chandra homomorphism is defined as the algebra homomorphism:
$
\gamma=\gamma(k) \:\Cal D \rightarrow
\Bbb C[E \otimes_{\Bbb R} \Bbb C]
$
by

$$
\gamma \: P = \sum_{\mu \ge 0} {e^{\mu(u)} p_{\mu}\big(\frac {\partial}
{\partial u_1},\ldots, \frac {\partial}
{\partial u_{n+1}}\big)}
\rightarrow \{\lambda \rightarrow p_{0}(\lambda +\rho)\}\;,
$$ where $\lambda=(\lambda_1, \lambda_2, \ldots , \lambda_{n+1})$ such
that $\lambda_1+\lambda_2+ \ldots + \lambda_{n+1} =0$ and $\rho$ is
defined as above.
\enddefinition

In particular,
$$
\gamma (L)(\lambda)= (\lambda +\rho, \lambda +\rho)- (2 \rho,\lambda
+\rho )=(\lambda, \lambda)-(\rho, \rho)
.$$
Also, one should notice, that
$$
\gamma (L)(w \lambda)= (w \lambda +\rho, w \lambda +\rho)- (2 \rho, w \lambda
+\rho )=(w \lambda,w \lambda)-(\rho, \rho)=(\lambda, \lambda)-(\rho,
\rho)
$$
\smallskip

\definition{0.5. Definition}(Action of the Weyl group $W$ on $\Cal
R$ and $\Cal D$)

The Weyl group acts on
$\Cal{R}$ by
$$
w \: {\frac{1}{1 - {e^{u_i-u_j}}}}
\to {\frac{1}{1 - {e^{u_{w(i)}-u_{w(j)}}}}}.
$$
One should notice that
$$
(1-e^{u_j-u_i})^{-1}= 1 - (1-e^{u_i-u_j})^{-1}
$$
and thus the action of the Weyl group is defined correctly.

The Weyl group acts on differential operators from $\Cal D$ as
$$
\multline
w \:
\sum_{i} f_{i}(u_1,\ldots, u_{n+1}) P_{i}\bigg({\frac{\partial}
{\partial{u_1}}},
\ldots,{\frac{\partial} {\partial {u_{n+1}}}}\bigg) \\
\rightarrow
\sum_{i} f_{i}(u_{w(1)},\ldots, u_{w({n+1})}) P_{i}\bigg({\frac{\partial}
{\partial{u_{w(1)}}}},
\ldots,{\frac{\partial} {\partial {u_{w({n+1})}}}}\bigg)
\endmultline
$$
\enddefinition
In particular, $L$ is Weyl group invariant differential operator.

\subhead {0.6. Symmetric  polynomials }\endsubhead
Consider symmetric polynomials in $u_1, \ldots , u_{n+1}$.
Let
$\sigma_1=\sum_{i=1}^{n+1} u_i \;,$
$\sigma_2=\sum_{i<j} u_i u_j \;,\ldots,$
$\sigma_{n+1}= \prod_{i=1}^{n+1} u_i $ be elementary
symmetric polynomials.

\proclaim{ 0.7. Theorem }
Any symmetric polynomial  $s =s(u_1, \ldots, u_{n+1})$ over $\Bbb C$
can be represented as a polynomial over  $\Bbb C$ of
elementary symmetric polynomials $\sigma_1, \ldots, \sigma_{n+1} \:$
$$s =s(\sigma_1, \ldots, \sigma_{n+1})$$
 and there is no algebraic relations
between elementary symmetric polynomials $\sigma_i$.
\endproclaim

\definition{ 0.8. Algebra of commuting differential operators }
Set  $\Bbb D=\Bbb D(k)$ for the algebra of all Weyl group invariant
differential operators in $\Cal D$ which commute with
operator $L=L(k)$:

$$\Bbb D=\{P \in \Cal D \;|\; [L,P]=0,\;
w(P) =P {\; \text{for any} \;} w \in W\}$$
\enddefinition
\smallskip
\remark{Remark}Restriction of Harish-Chandra homomorphism $\gamma$
from
$\Cal D$ to $\Bbb D$ will be also called Harish-Chandra homomorphism.
\endremark

\smallskip
\proclaim{0.9. Theorem} 1). Differential operator
 $P = \underset{\mu \ge0}\to{\sum} e^{\mu(u)}p_{\mu}(\frac {\partial}{\partial
u})$ in $\Cal D$
commutes with $L\: [ L, P] = L \circ P -P\circ L =0$ if and only if polynomials
$p_{\mu}(\lambda)$ satisfy the recurrence relations:
$$
(2 \lambda -2 \rho +\mu,\mu) p_{\mu}(\lambda) =
2 k\sum_{\alpha \in R_{+}} \sum_{j=1}^{\infty} \{(\lambda +\mu -j
\alpha, \alpha)p_{\mu -j \alpha}(\lambda)  - (\lambda, \alpha) p_{\mu -j
\alpha}(\lambda + j
\alpha)\} .
$$

2). For any symmetric polynomial $\sigma =\sigma(\lambda)$ set
$p_{0}(\lambda) = \sigma (\lambda - \rho)$. Then all
the recurrence relations above can be solved and obtained differential
operator is in $\Cal D$,  Weyl group invariant and thus in $\Bbb D$.
Furthermore,   all elements of $\Bbb D$ are obtained in this way, or in
other words, Harish-Chandra homomorphism
$
\gamma=\gamma(k) \:\Bbb  D \rightarrow
\Bbb C[E \otimes_{\Bbb R} \Bbb C]^W
$ is surjective.

3). All differential operators constructed in 2). commute not only with
$L$, but as well  with each other.
\endproclaim

For example, to 1 we assign in this way trivial differential operator
, namely, identity, to
$\sigma_1=\sum \lambda_i$ we assign $\sum \frac {\partial}{\partial
u_i}$, to $\sum \lambda_i^2$ we assign $L+(\rho, \rho)$ , and so on.

\remark {Remark 0.10 } For explicit form of generators of $\Bbb D(k)$ reader
should consult with [ 38,39,40].
\endremark

\subhead { Hypergeometric system of differential equations }\endsubhead
\definition{0.11. Definition} ( Hypergeometric system of differential
equations)
 The system of differential equations

$$
P \phi = \gamma (P)(\lambda) \phi
\quad {P} \in \Bbb D, \quad \lambda =(\lambda_1,\ldots,\lambda_{n+1})\;
\text{s.t.}\; \lambda_1 +\ldots + \lambda_{n+1}=0
$$

is called the system of hypergeometric (partial)
differential equations.
\enddefinition
One may consider this system in variables $u_i$ and in variables $z_i$.
To obtain a system in $z_i$ one should replace $\frac
{\partial}{\partial u_i}$ by $ z_i \frac
{\partial}{\partial z_i}$ and $(1-e^{u_i-u_j})^{-1}$ by
$\frac{z_j}{z_j- z_i}$ correspondingly.

\definition{ 0.12. Definition }Set  $H^{reg}$ as a (n+1)-dimensional complex
plane $\Bbb C^{n+1}=\{(z_1,\dots,z_{n+1})\}$ with  diagonals
$\{z_i=z_j\},\;i<j$, and with coordinate
hyperlanes
$\{z_i=0\}$ being deleted:
$$
H^{reg}=\Bbb C^{n+1}\setminus \bigg(\bigcup \{z_i=z_j\}  \bigcup
\{z_i=0\}\bigg)
$$
\enddefinition

\remark{Remark} This configuration is also called an affine configuration. It
became classical now and was studied by Arnold,Brieskorn, Orlik,
Solomon, Schechtman, Varchenko,...
\endremark

\proclaim {0.13. Theorem } (Holonomicity on $H^{reg}$)
Locally on $H^{reg}$  solution
space of the
system of hypergeometric equations has dimension equal to the order
$|W|$ of the Weyl group $W$ and consists of analytic functions.
\endproclaim

\smallskip
\subhead{0.14. Asymptotic solutions}\endsubhead
Operator $L$ permits the following asymptotic  expansion :

$$
L=
\sum {\frac {\partial^2}{\partial u_i^2} } -2\sum_i (\rho,
e_i) \frac{\partial}{\partial u_i}-
2\sum_{i < j } k \sum_{m=1}^{\infty}  e^{m(u_i-u_j)}
\big(\frac {\partial}{\partial u_i}-\frac {\partial}{\partial u_j}\big).
$$

Consider  solutions of
$$
L \phi =((\lambda , \lambda) -(\rho, \rho)) \phi
$$

of the form

$$
\phi= \phi(\mu,k,e^{u})= \sum_{\nu \ge \mu} \Gamma_{\nu}(\mu,k) e^{\nu(u)}
\quad ,
$$

where  $\mu = w \lambda +\rho$.

Then
coefficients $\Gamma_{\nu}(\mu,k)$ satisfy to
 Freudenthal-type recurrence relations
$$
\{(\nu -\rho, \nu -\rho)-(\mu -\rho, \mu -\rho) \}
\Gamma_{\nu}(\mu,k)=
 2 \sum_{\alpha \in R_{+} }k \sum_{j=1}^{\infty} (\nu -j \alpha, \alpha)
\Gamma_{\nu- j \alpha}(\mu,k).
$$

For generic $\lambda$ $((\lambda, \alpha^{\vee}) \notin \Bbb Z, \;
\alpha \in R)$ solutions $\phi(w \lambda +\rho,k,e^u)$ are
linearly independent and provide a basis of linear space of solutions
to the whole system of hypergeometric equations.

These solutions  $\phi(w \lambda +\rho,k,e^{u})$  (or $\phi(w \lambda
+\rho,k,z)$ )
are called {\bf{ asymptotic
solutions}}.
\smallskip
{\it {So instead of dealing with the whole hypergeometric system, one
can deal with differential operator of second order and asymptotic
solutions only.}}

\subhead {0.15. Integral representation for  the solutions to the
hypergeometric system of
differential equations} \endsubhead

A.Matsuo [15] proved that
the hypergeometric system of
differential equations in the  case of root system of type $A_n$
( which is only being considered in this paper)
is related to the Knizhnik-Zamolodchikov
equation in conformal field theory, cf. also [16 ].
In particular this implies that solutions of the system are
represented as certain multidimensional integrals cf.  [ 1 ].

Namely, take integration variables indexed as
$\{t^{(j)}_l;\; j=1,\ldots,n,\;
l=j,\ldots,n\}$. Take a total ordering $<$ in the set of indices $I$.

$$
\align
\omega(z,t)= &\prod_{(j,l) \in I} e^{-{k \over 2}(\lambda, \alpha_j)
t_l^{(j)}}
\prod_{1 \le i \le {n+1}}e^{{k \over 2}{(\lambda, \Lambda_1)}z_i}\\
&\times \prod_{{1 \le i \le n+1}\; {1 \le l \le n}} (\sinh \frac{z_i
-t^{(1)}_l}{2})^{-k} \prod_{(j,l) < (j^{\prime}, l^{\prime})} (
\sinh \frac{t_l^{(j)}-
t^{(j^{\prime})}_{l^{\prime}}}{2})
^{k(\alpha_j, \alpha_j^{\prime})}
\endalign
$$
Then
$$\int_{\Gamma} \omega(z,t) \phi(z,t) dt_1^{(1)} \ldots dt_{n}^{(n)}$$
where $\phi(z,t)$ is a polynomial of $\coth \frac{t_l^{(j)}-z_i}{2}$
and $\coth \frac{t_l^{(j)}-t_{l^{\prime}}^{(j^{\prime})}}{2}$ written in a
complicated manner in [ 1] for a certain choice  of $\Gamma$
provide solutions for hypergeometric system of differential equations.
Contours for asymptotic solutions and for zonal spherical function
itself are not described in [15 ].

In this paper we provide another integral representation for the
solutions cf. theorem 6.3. Namely, we use
$$
\align
\omega(z,t)= & {\prod (z_i -z_j)^{1-2k}}\\
&\times \prod_{(j,l) \in I} {t_l^{(j)}}^{(\lambda +\rho,- \alpha_j)}
\prod_{1 \le i \le {n+1}}{z_i}^{{(\lambda +\rho, e_1)}}\\
&\times \prod_{{1 \le i \le n+1}\; {1 \le l \le n}}
 ( {z_i-t^{(1)}_l})^{k-1} \prod_{(j,l) < (j^{\prime}, l^{\prime})}
({t_l^{(j)}-
t^{(j^{\prime})}_{l^{\prime}}})
^{(1- k)(\alpha_j, \alpha_j^{\prime})}
\endalign
$$
Let's imbed $V$ into $ \Bbb R \times V$ as $v \to (0,v)$. Let
$\alpha_0= e_0-e_1$.
Then our form $\omega(z,t)$ can be written as:
$$
\align
\omega(z,t)= & \prod_{i<j}(z_i -z_j)^{-1}
 {\prod (z_i -z_j)^{(1-k)(\alpha_0, \alpha_0)}}\\
&\times \prod_{(j,l) \in I} {t_l^{(j)}}^{(\lambda +\rho,- \alpha_j)}
\prod_{1 \le i \le {n+1}}{z_i}^{{(\lambda +\rho, -\alpha_0)}}\\
&\times \prod_{{1 \le i \le n+1}\; {1 \le l \le n}}
 ( {z_i-t^{(1)}_l})^{(1-k)(\alpha_0,\alpha_1)} \prod_{(j,l) < (j^{\prime},
l^{\prime})}
({t_l^{(j)}-
t^{(j^{\prime})}_{l^{\prime}}})
^{(1- k)(\alpha_j, \alpha_j^{\prime})}
\endalign
$$

Solutions to the hypergeometric system are given by the integral

$$
\int \omega(z,t) dt_1^{(1)} \ldots dt_{n}^{(n)}
$$
over appropriate contour of integration.
Here our $z_i$ is $e^{z_i}$ of Matsuo.
\smallskip
In the paper we use another indexation of variables of integration
$\{t_{i,j}; i=1,\ldots,j\;, j=1,\ldots,n \}$. Variables $t_{i,j}$
have a nice geometric origin  in elliptic coordinates cf.[25], also
they are used in [25] in the Plancherel theorem.

\remark {Remark 0.16 }
We obtained the form in the
following way.
 Zonal spherical function is defined as $\phi(g)= (T_g \xi,\xi)$,
where $T_g$ is a unitary representation of $g \in G$,  $\xi$ is
 invariant vector of maximal compact subgroup cf.[25].
 The zonal spherical function for $SL(n,\Bbb C)$ was calculated in
[25] and using the same methods for $SL(n,\Bbb R)$ in [35]. Essentially,
the so-called elliptic coordinates are used. The  case of  $SL(n,\Bbb C)$
corresponds to $k=1$ in [18] and case of $SL(n,\Bbb R)$ corresponds to
$k=\frac{1}{2}$. So if one considers integral representations for zonal
spherical functions for $k=1$ and for $k=\frac{1}{2}$ (one should also
replace ${\delta_i}^2$ by $z_i$ and $i \lambda \over 2$ by $\lambda$, i.e.
adopt normalizations of [18] )
 and then connects
powers of factors linearly on $k$, one obtains exactly the above formulas.
Also, in calculation of zonal spherical function integration is taken
over a distinguished cycle. Its role  is discussed  in [33].
\endremark
\smallskip

\subhead{0.17.  Organization of the paper} \endsubhead

We describe contours  for  integration for asymptotic solutions
$\phi(w \lambda +\rho, k,z)$ and
 enumerate them by the elements
of symmetric group. Namely, for each element $w \in S_n$ we put into
correspondence a diagram, cf. sections 1,2.

For the diagrams we borrow very
convenient graphical notation of [2 ].
We also discuss the interrelation of diagrams and Gelfand-Zetlin
patterns in section 3 (Weyl group orbit of the lowest weight is
described with the help of diagrams).

The system of contours $\Delta_w= \Delta_w(z)$ for integration for each diagram
( and thus for
each element $w \in S_n$) is described in section 4, cf. definition 4.3.

Section 5 is devoted to description of the choice of the phase of
multivalued form $\omega$ over cycle $\Delta_w$.

In section 6 we  calculate the leading asymptotics
 $ (w\lambda +\rho)$ of the integral  and leading
coefficient, and
prove that integrals satisfy to second order differential equation
 $L \phi =((\lambda , \lambda) -(\rho, \rho)) \phi$
and thus to the whole system of hypergeometric equations.
Finally (theorem 6.8) , we
calculate the value of the integral over described cycles when all the
arguments collapse to the unity. We use result of [21 ], and it is
consistent with [17 , 33 ].
\remark {Remark 0.18 }
In the case of $k=1$ our integral simplifies dramatically ( boils
down to integrals considered in [25]), namely
only affine part (or Mellin part) of the integral survives. Also,
the zonal spherical function in this case  is equal to the quotient of
some determinant by $\prod (z_i -z_j)$. The determinant is a sum
of $n!$ terms with the signs plus or minus. So one can try to select
each term from the determinant, i.e. to find an appropriate domain of
integration. This can be easily done. Modification of these domains to
the case of generic $k$
is done in the definition 4.3.
\endremark

\remark {Remark 0.19}
Our integral representation for the solutions $\int \omega(z,t) dt$
is in well agreement with the formula of Harish-Chandra for zonal
spherical functions:
$$\phi_{\lambda}(g)= \int e^{(i \lambda +\rho) (A(kg))}dk, \quad g \in G$$
For more details about formula of Harish-Chandra and application of
the stationary phase method to the formula cf.  [20] and  [43],
correspondingly.

\endremark

\head 1. Diagrams\endhead

 The notion of a diagram was introduced in [1]
in the context of Knizhnik-Zamolodchikov equation, and later similar
notion was introduced in [2] in the context of
 multidimensional determinants. We will use the notion of a diagram in the
form of [2], in particular, we borrow the graphical notation of [2].

 Compare also the  diagrams with  sequences Seq of  [3].
\medskip
Fix some positive integer $n$. Consider the set of $\frac{n(n+1)}2$
points, indexed by pairs of integers
 $\{(i,j)|$ $i=1, \ldots,j$, $j=1,\ldots, n \}$.
It is helpful to organize the points in the form
of a pattern, so that points are divided in $n$ rows, $j$th row is formed
by points $\{(i,j)| \; i=1, \ldots,j\} $; point $(i,j)$ is located under and
between points $(i,j+1)$ and $(i+1, j+1)$ (fig. 1a).

\midinsert\vskip 4cm
\includegraphics{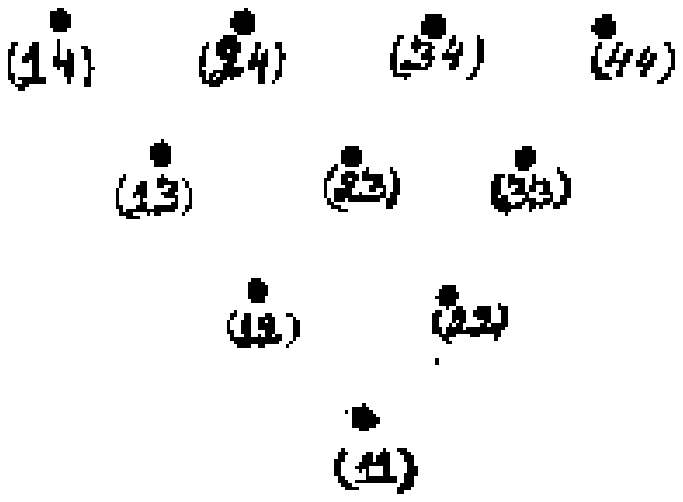}
\botcaption{Figure 1a}
$n=4$.\quad Points $ (i,j)$ organized in a pattern;
$j$ is the number of the row, $i$ is the number in the row
\endcaption
\endinsert

Now mark with a cross one point in each row. Let $\{(i_j,j)\}$ be the
subset of marked points (fig. 1b).

\midinsert\vskip 4cm
\includegraphics{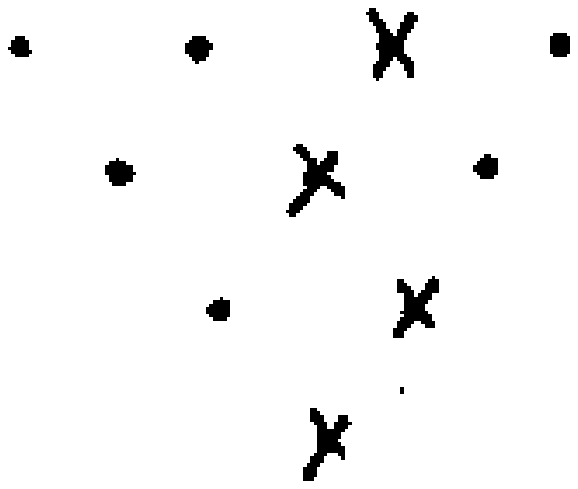}
\botcaption{Figure 1b}
One point is marked in each row.
\endcaption
\endinsert

Finally, draw an arrow for each point $(i,j)$ with the source in this
point $(i,j)$ and target $tar(i,j)$ in the next $j+1$th row defined as:
$$
tar(i,j)=
\cases
(i,j+1), &\text{if }i<i_{j+1}\\
(i+1,j+1), &\text{if }i \geq i_{j+1}\text{\footnotemark}
\endcases
$$
\footnotetext{Recall that $\{(i_j,\;j) \} $ is the set of marked points and
 $ (i_{j+1},\; j+1)$
is the only marked point in $j+1$th row.}

Note: neither arrow has a marked point as its target.

In this way one obtains fig 1c.

\midinsert\vskip 4cm
\includegraphics{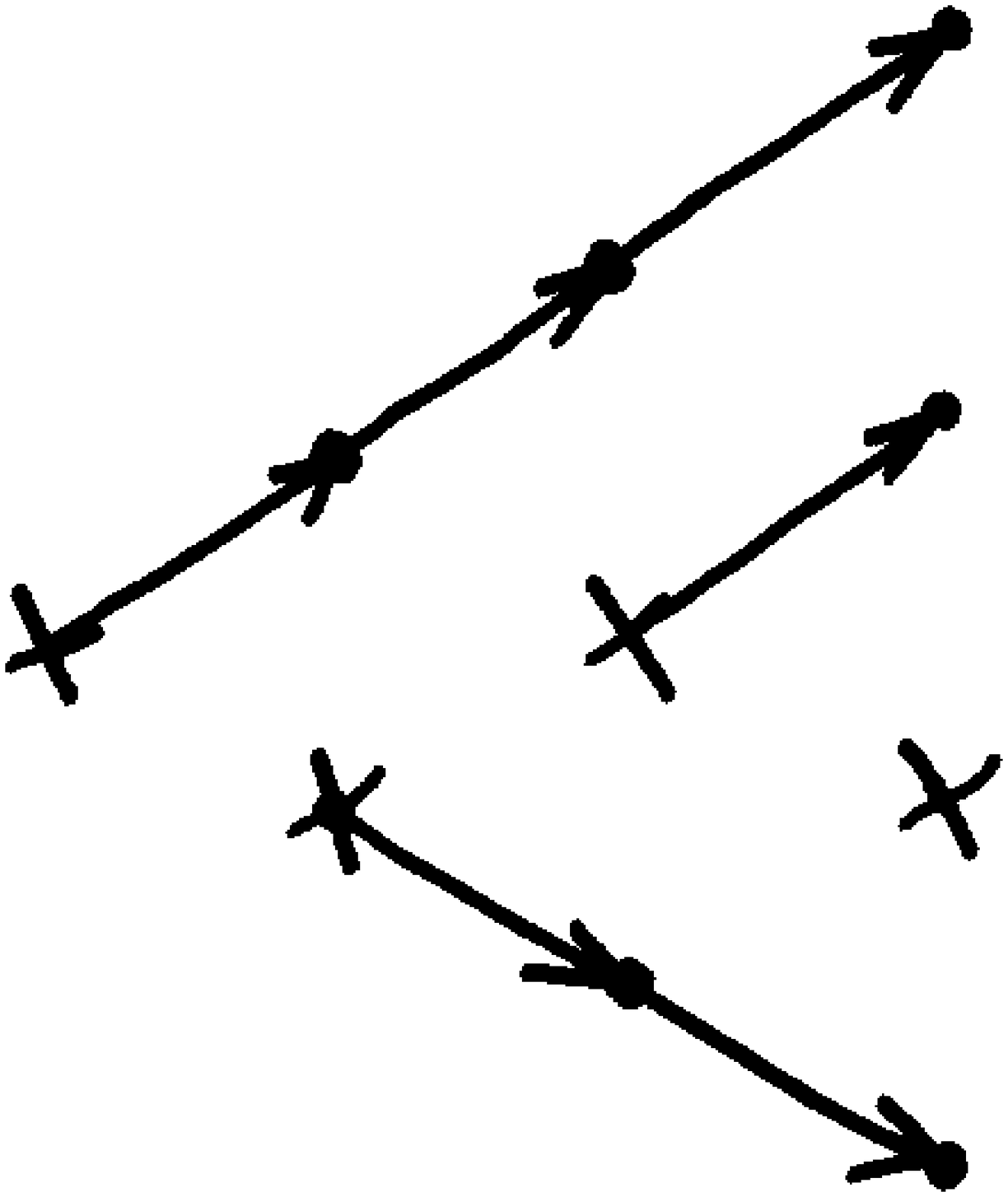}
\botcaption{Figure 1c}
Example of a diagram.
\endcaption
\endinsert

\definition{Definition 1.1}
A triple $$ (\{i,\;j\},\quad \{i_j,\;j\},\quad tar) $$
 consisting of:

set of points
 $\{(i,\;j)| \quad i=1,\ldots,j\; ; \; j=1, \ldots,n\} $,

set of marked points $\{ (i_j,\;j) \; | \quad j=1,\ldots,n \} $

and function
$\quad tar \quad $ defined above will be called {\bf a diagram.} \enddefinition

\remark {Remark 1.2} One can see that a diagram is determined by the set of
marked points.
\endremark

One easily notices that there are exactly $n!$ different diagrams with $n$
rows. So it is desirable to enumerate diagrams by elements of
Symmetric group $S_n$. Calculations connected with the integral led us
to the following propositions and theorems. One may consider them as a
peculiar show up of quantum groups.

\head{2. Enumeration of diagrams by elements of Symmetric group}\endhead

\definition{Definition 2.1} Consider a diagram. For a point $(i,\;j) $,
provided
$(i,j)$ is not a marked point, $j=2, \ldots,n $, define its {\bf source} as
 $(i_1,j-1)$,s.t. $tar(i_1, j-1)=(i,j)$, i.e.
 $\quad source =tar^{-1} \quad $
, which
 is
defined on the set of nonmarked points.
\enddefinition

Consider a diagram with $n$ rows and a point $(i,n)$ in the $n$th row.
Consider also along with $(i,n)$ its
 $\op{source}(i,n)$, $\op{source}(\op{source}(i,n))$,
and so on until we get a marked point. Say,$\quad source^{w(i)-1}(i,n)$
is a marked point.

\proclaim{Proposition 2.2} Correspondence $i\rightarrow w(i)$ correctly defines
 permutation of
numbers $1,\ldots,n$, i.e. an element of symmetric group $S_n$.
This is a one-to-one correspondence between diagrams with $n$ rows and
symmetric group $S_n$.
\endproclaim

\remark{Remark 2.3} If $(i_n,\;n)$ is marked point in $n$th row of a
 diagram with $n$
rows, then $ w(i_n)=1$. Also if $ (i_{n-1},n-1)$ is a marked point in
$ n-1$th row
, then $ w(i_{n-1})=2 $, provided $i_{n-1} < i_n$ and $ w(i_{n-1}+1)=2$,
 provided $i_{n-1} \ge i_n$
and so on.
\endremark

\remark {Remark 2.4} In other words one can describe the correspondence between
diagrams and elements of Symmetric group as follows.
Consider a diagram as an oriented graph and forget orientation.
For $i=1, \ldots,n$ define $ w(i)$ as the number of points in the connected
component of the point $(i,n)$.
\endremark

Symmetric group $S_n$ has standard generators
 $\sigma_1,\;\sigma_2,\ldots, \sigma_{n-1}$, where
 $\sigma_i$ permutes $i $ and $ i+1$.

\definition{Definition 2.5} The length $l(w)$ of an element $w \in{ S_n}$ is
the minimal
integer $p\ge 0$, s.t. $w$ admits a presentation

$$ w = \sigma_{i_1}\sigma_{i_2} \ldots \sigma_{i_p}  \;.$$
Any presentation of $w$ as a product of $p=l(w)$ generators is called
a reduced presentation.

\enddefinition

\proclaim {Theorem 2.6}
 Let a diagram
$$
(\{(i,j)| \quad i=1, \ldots,j,\; j=1, \ldots,n\}, \quad \{(i_j,j) | \quad
j=1, \ldots,n \}
 \quad, \; tar)
$$
corresponds to an element $w\in{ S_n}$. Then the length $ l(w)$ of an
element $w$ is equal to:
$$
l(w)=\sum_{j=1}^n{(i_j -1)}
$$
\endproclaim

\remark{Remark 2.7} Arrow with the source $t_{ij}$ has either  $t_{i,j+1}$
or  $t_{i+1,j+1}$ as its target. So one can say that the arrow is {\bf{to
the left}} or {\bf{to the right}}. Hence the theorem 2.6 says that $l(w)$ is
equal
to the number of arrows , which are to the left.
\endremark

\smallskip
The following nice theorem is well known, cf. [4].
\smallskip

\proclaim{Theorem 2.8} For symmetric group $S_n$ the following identity holds:

$$
\sum_{w\in{S_n}} {q^{l(w)}} =
 {(1-q)\over(1-q)}{(1-q^2)\over(1-q)}\ldots{(1-q^n)\over(1-q)} \quad .
$$
\endproclaim

\smallskip
Though throughout this paper we mainly using Symmetric group $S_n$
for purposes of dealing with diagrams, for a moment we would like to
use the diagrams to obtain an extension of the above theorem.

\proclaim{Theorem 2.9}(Multiparametric deformation of symmetric group) Denote
by ${\{i_j(w),j\}}$ set of marked points
of
 diagram
corresponding to $w\in{ S_n}$. Then the following identity holds:
$$ \sum_{w\in{S_n}}\prod_{j=1}^n {q_j^{i_j(w)-1}} =
 {(1-q_1)\over(1-q_1)}{(1-q_2^2)\over(1-q_2)} \ldots
{(1-q_n^n)\over(1-q_n)}.
$$
\endproclaim

\smallskip
This theorem is an easy application of the formula for the sum of
geometric progression
$$
1+x+x^2 +x^3 + \ldots + x^{n-1} = {{1- x^n}\over{1-x}}.
$$

\proclaim{Theorem 2.10}
Suppose $w$ corresponds to a diagram with set of marked points
$\{(i_j,j)\}$. Then $w$ admits a presentation:

$$w = w_n w_{n-1} \ldots w_2 w_1 \;, $$
where

$$
w_k =\cases
\sigma_k \sigma_{k+1} \ldots  \sigma_{i_{n-k+1} +k-2} , &\text{if }i_{n-k+1} >
1\\
id , &\text{if }i_{n-k+1}=1
\endcases
$$
 Moreover, presentation of $w$ in terms of generators
$\sigma_i, \quad i=1,\ldots , {n-1}$, naturally obtained from the above
presentation by expanding $w_k$ in terms of generators
(omitting first those $w_k$ which are identities)  is a reduced
presentation of $w$,
in particular, it contains exactly $\sum_{j=1}^n {(i_j-1)} =l(w)$
factors.
\endproclaim

\subhead {2.11. Proof of theorem 2.6} \endsubhead Theorem 2.10 provides us for
each
 $w \in S_n$ a presentation with $\sum_{j=1}^n {(i_j(w)-1)}$ factors.
Consequently, for each $ w$  $\sum_{j=1}^n {(i_j(w)-1)} \le l(w)$.
But comparison of theorems 2.8 and 2.9 shows that all inequalities are
in fact equalities. Thus theorem 2.6 is proved.
\smallskip

\definition{ Definition 2.12 }
(Partial ordering on Symmetric Group)
Consider two elements of Symmetric Group $S_n$, say, $w_1$ and $w_2$.
Let diagram corresponding to $w_1 \; (w_2)$ has $\{i_j(w_1),j\}$
($\{i_j(w_2),j\}$)  as the set of marked points.
We say that {\bf{ $w_1$ is less than or equal to  $w_2$}}:
$w_1 \le w_2$
if for each $j=1,\ldots,n$
 $$i_j(w_1) \le i_j(w_2).$$
\enddefinition
\smallskip

\proclaim{ Proposition 2.13 } \endproclaim
a). For each element $w$ of Symmetric group
$S_n$ (
$w \in S_n$)

there are exactly
$$
\prod_{j=1}^n {(j-i_j(w)+1)}
$$
elements which are greater than or equal to $w$.
\smallskip
b).
$$
\sum_{w^{ \prime} \ge w}
q^{l(w^{\prime})}=
q^{l(w)}\prod_{j=1}^n
{
\frac
{1-q^{j-i_j(w)+1}}
{1-q}}
$$
\smallskip
c). For each element $w$ of Symmetric group
$S_n$ (
$w \in S_n$)
there are exactly
$$
\prod_{j=1}^n {i_j(w)}
$$
elements which are less than or equal to $w$.
\smallskip
d).
$$
\sum_{w^{ \prime} \le w}
q^{l(w^{\prime})}=
\prod_{j=1}^n
{
\frac
{1-q^{i_j(w)}}
{1-q}}
$$
\smallskip
e). If $w \le w^{\prime}$, then $l(w) \le l(w^{\prime})$.
\medskip

\head 3. Gelfand-Zetlin patterns \endhead

In  [5] finite-dimensional
representation of $gl(n)$ with highest weight
 $(m_n,m_{n-1},\ldots, m_1)$ s.t. $m_n \ge m_{n-1} \ge \ldots \ge m_1$
is proved to have a nice basis , which elements are enumerated by the so-called
Gelfand-Zetlin patterns, i.e. set of numbers
 ${m_{pq}, \; p\le q ,\;q=1,\ldots, n-1}$
arranged in the following pattern:

\midinsert
$$
\matrix
m_{1n}&&m_{2n}&&\ldots&&\ldots&&m_{nn}\\\\
&m_{1,n-1}&&m_{2,n-1}&&\ldots&&m_{n-1,n-1}&\\\\
&&\ldots&&\ldots&&\ldots&&\\\\
&&&m_{1,2}&&m_{2,2}&&&\\\\
&&&&m_{1,1}&&&&
\endmatrix
$$
\botcaption{Figure 2}
Gelfand-Zetlin pattern
\endcaption
\endinsert

The numbers are arbitrary integers which satisfy inequalities
$m_{p,q+1} \le m_{pq} \le m_{p+1,q+1}\;, \quad p=1,\ldots,q,
\quad q=1,\ldots, n-1 $ . We changed the usual
inequalities to the opposite. The numbers
 $ m_1, \ldots,m_i,\dots,m_n$
 which define the representation are denoted by $m_{in}$
and placed in the nth row.
\smallskip
Consider a diagram $\{(i,j),\quad (i_j, j) ,\quad tar \}$
with n rows, which corresponds to an element $w\in{S_n}$.

For the highest vector $m_n \ge m_{n-1} \ge \ldots \ge m_1$
and  $w\in{S_n}$ we put into correspondence a Gelfand-Zetlin
pattern, which is uniquely defined by the relations :

$$m_{in}=m_n, \quad m_{pq} = m_{tar(p,q)}.$$

\proclaim{Theorem 3.1} The above correspondence
$w\rightarrow \{m_{pq}(w)\}$ is actually the action of
$w\in{S_n}$  on the lowest weight $(m_1,m_2,\ldots, m_n)$ of
representation of $gl(n)$, i.e. $\{m_{pq}(w)\}$ is the only vector
of weight $(m_{w(1)}, m_{w(2)}, \ldots, m_{w(n)}).$
\endproclaim

\remark{Remark 3.2} This section is aimed to emphasize the
philosophical relation between Harish-Chandra decomposition and
BGG resolution. \endremark

\smallskip

\head 4. Cycles for asymptotic solutions\endhead

\definition{Definition 4.1}
 By a  {\bf bump function} we mean a $C^1$ function
$$
f_{\epsilon}\: [0,1]\rightarrow[0, \epsilon],
$$
 such that
$$ f_{\epsilon}(x)=0 \quad \Longleftrightarrow \; x=0 \;or \; x=1
\quad .$$
\enddefinition

\midinsert\vskip 5.5cm
\includegraphics{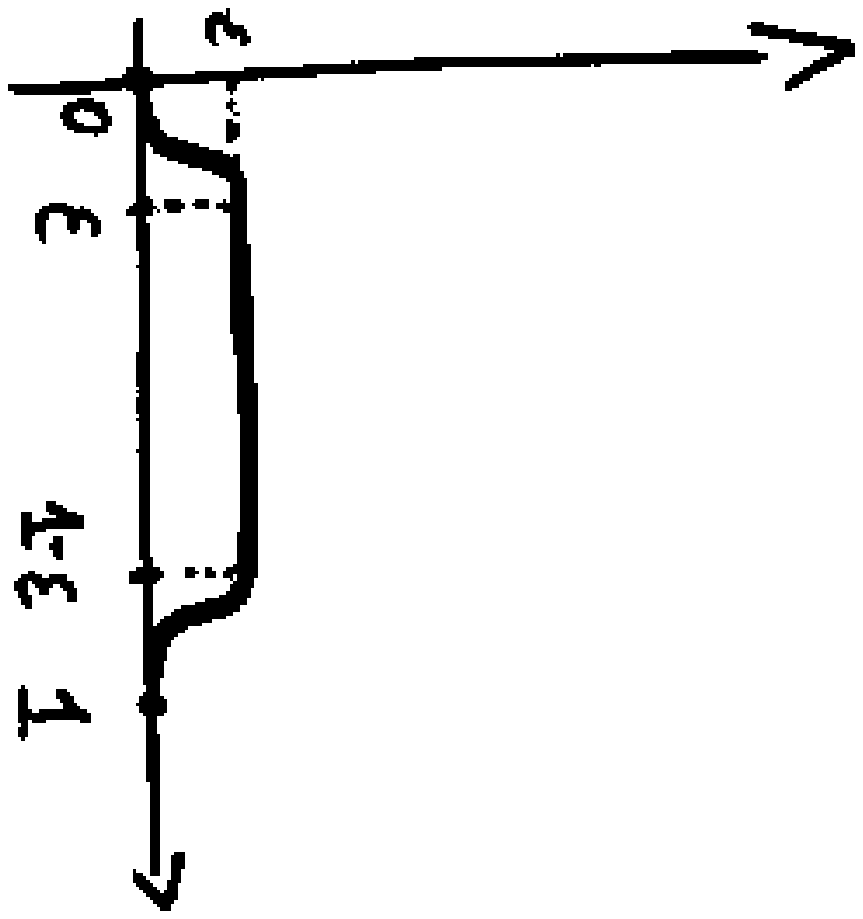}
\botcaption{Figure 3}
Bump function
\endcaption
\endinsert

\remark{Remark 4.2}We assume $\epsilon$ to be sufficiently small as needed.
\endremark

\medskip
Consider a diagram

 $$(\{(i,j)|\quad i=1,\ldots,j, \;\;j=1, \ldots,n+1\},
 \quad \{(i_j,j)|\quad j=1, \ldots,n+1 \},
\quad tar)$$
 corresponding to an element $w\in{S_{n+1}}$.

With each point $(i,\;j)$ one associates a variable $t_{i,j}$.
Later variables
 $$\{t_{i,\;j}|\; i=1, \ldots,j;\; j=1, \ldots,n \}$$
 will be variables
of integration, while variables
 $$t_{i,n+1}, \quad i=1, \ldots, n+1$$
will have the meaning of arguments. So let
$$ \{t_{1,n+1}, \;t_{2,n+1}, \; \ldots,\; t_{n+1,n+1} \}$$
be fixed and
$$0 < | t_{1, n+1}| < |t_{2, n+1}| <\ldots < |t_{n+1, n+1}| \; .$$
\smallskip

\remark{NOTATION} Variables $ t_{1, n+1}, \ldots,t_{n+1, n+1}$ will be also
denoted by
$z_1,\ldots, z_{n+1}$, correspondingly.
\endremark

\definition { Definition 4.3}
For each
 $ \quad t_{i, \; j}, \quad i=1, \ldots,j;\quad j=1, \ldots,n$ define a contour
$\quad t_{i,\;j}(\tau_{i,\;j})$, $ \quad \tau_{i,\;j}\in[0,\;1] \quad $
as
$$
t_{i,j}(\tau_{i,j})=e^{2\pi i\tau_{i,j}}(1-f_{\epsilon}(\tau_{i,j}))
\; t_{tar(i,j)} \; .
$$
\enddefinition
See fig. 4a, 4b, 5a, 5b, 5c, 5d, 5e, 5f.

\smallskip
 The described system of contours corresponding to $w\in{S_{n+1}}$
will be denoted by $\Delta_w=\Delta_w(z)$ and used as a cycle for integration
of
the form described in the next section.
Compare also with [6,7,8,9,10,11,12,13,32].

\remark{Remark 4.4} In $\Delta_w$  one has $|t_{ij}| \le
|t_{tar(ij)}|$, cf. the asymptotic zone of [14].
\endremark

\head 5. Multivalued form $\omega_w$\endhead

\smallskip
Let
$$
(\{(i,j)|\quad i=1, \ldots ,j,\; j=1, \ldots,n+1\}, \quad \{(i_j,\;j) | \quad
j=1,...,n+1 \},\quad tar)$$
 be a diagram corresponding to an element $ w\in{S_{n+1}}$.

\definition{Definition 5.1}
Define $$x \: \Bbb Z\times \Bbb Z\rightarrow{\Bbb Z}$$ to be the projection on
the
first factor, i.e.

$x(a,b)=a$.
\enddefinition

\definition{Definition 5.2} Define $ \omega_w $ as
$$
\align
\omega_w:=&
 \prod_{i=1}^{n+1} t_{i,n+1}^{\lambda_1 +{k n\over 2}}
 \prod_{i_1 \ge i_2}(t_{i_1,n+1}-t_{i_2,n+1})^{1-2k}
 \prod_{j=1}^n \{
 \prod_{i_1 < x(tar(i,j))} ( t_{ij} - t_{i_1,j+1} )^{k-1}\\
 &\times\prod_{i_1 \ge x(tar(i,j))} { ( t_{i_1,j+1}- t_{ij} )^{k-1}}
 \prod_{i_1>i_2} {(t_{i_1,j}-t_{i_2,j})^{2-2k} }\\
 &\times\prod_{i=1}^{j} {t_{ij}^{\lambda_{n-j+2}-
  \lambda_{n-j+1} -k} } \}\quad
 {dt_{11} dt_{12} dt_{22} \ldots dt_{nn} }
\endalign
$$
cf. [1, 9, 10, 15, 17].
\enddefinition

In $\Delta_w$ we have $t_{ij}=t_{ij}(\tau_{ij}) ,\quad i=1, \ldots ,j,
\quad  j= 1, \ldots, n $($ t_{i,n+1} \quad i=1,\ldots,n+1 $are fixed)
and the phase of factors
in the formula of $\omega_w$ should be chosen so that it goes to zero
as $\tau_{ij}$ approach to zero, provided $k$ and
 $\lambda_1,\; \lambda_2,\; \lambda_3,\ldots,\; \lambda_{n+1}$ are
real, cf. Remark 4.4.

\remark {Remark 5.3} Index $w$ of the multivalued form above is, of course, a
convention about the phase of the same form over the cycle $\Delta_w$ \; .
\endremark

\head 6. Main theorems \endhead

\proclaim{ Theorem 6.1} Let $w\in{S_{n+1}} $. Then the integral of the
multivalued
form $\omega_w$ over cycle $\Delta_w$ gives an asymptotic solution

$\phi(w \lambda+\rho, k,z) \:$

$$
\phi(w \lambda+\rho, k,z)=
\int_{\Delta_w(z)} \omega_w=
 a(w) z^{w \lambda +\rho}(1+ \ldots )
$$

where

$$
 z^{w{\lambda +\rho}}= { z_{1}^{ { \lambda_{w(1)} } +{k
n\over2}} \;\;
 z_{2}^{ {\lambda_{w(2)} } + {k(n-2)\over2}} \; \ldots
 z_{n+1}^{{\lambda_{w(n+1)}
 -{ kn\over2}}}
 }
$$
and
$$
\align
a(w)=\prod_{\alpha\in R_{+}}&
\frac{\Gamma((-w{\lambda}, {\alpha^\vee})) \sin( \pi(-w \lambda, \alpha^\vee))}
{\Gamma((-w{\lambda},{\alpha^\vee)}+k)}\\
&\times e^{-2\pi i(\lambda, \delta)} e^{-\pi i(k-1)l(w)}
 \Gamma(k)^{\frac{n(n+1)}2} (2i)^{\frac{n(n+1)}2} \quad .
\endalign
$$
\endproclaim

The theorem is proved by induction on number of rows of a diagram.

\vskip 1cm

\subhead{Mechanism of induction 6.2}\endsubhead
The mechanism of induction is based on the following simple observation.
Let a diagram
$$
 (\{i,\;j\},\quad \{i_j,\;j\},\quad tar| \; i=1,\ldots,j, \; j=1,\ldots,n+1)
$$
with n+1 rows corresponds to an element $w \in S_{n+1}$.

Consider a diagram with $n$ rows, which is obtained from the diagram
with  $n+1$ rows by deleting the $n+1$ row :
$$
 (\{i,\;j\},\quad \{i_j,\;j\},\quad tar| \; i=1,\ldots,j, \; j=1,\ldots,n)
$$ and suppose that it corresponds to $ w^{\prime} \in S_n$.
Then one has:
$$
w(i)=
\cases
 w^{\prime} (i)+1, &\text{if }i<i_{n+1}\\
 w^{\prime}(i -1) +1, &\text{if }i >i_{n+1}\\
 1,  &\text{if }i =i_{n+1}
\endcases
$$
cf.  also theorems 2.10 and 2.6 and fig. 6.

\midinsert\vskip 5cm
\includegraphics{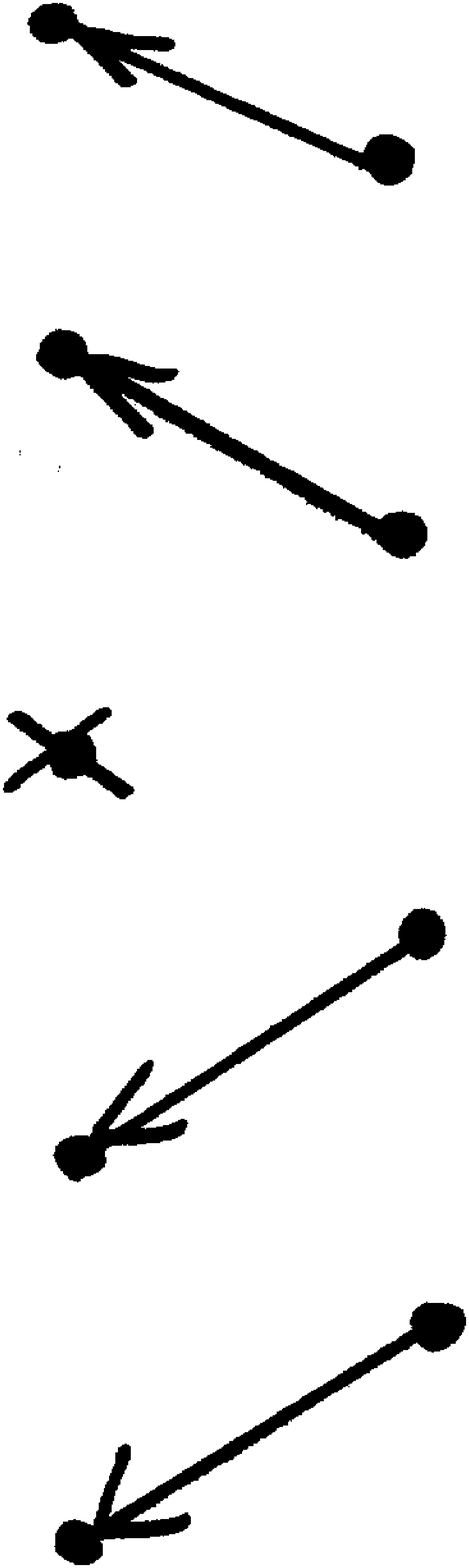}
\botcaption{Figure 6}
Illustration of induction
\endcaption
\endinsert

\smallskip

\proclaim{ Theorem 6.3} The integrals of multivalued forms  $\omega_w$ over
cycles $\Delta_w$
satisfy to the following differential operator of second order:

$$
\align
\{\sum_{i=1}^{n+1} { ({{z_{i}{\partial\over{ \partial z_{i} }
}}})^2}  - k \sum_{i <j} &{
{\frac{z_{j}+z_{i}}{z_{j}-z_{i}}}
 (z_{i} {\partial \over
\partial z_{i}}- z_{j} {\partial \over \partial z_{j}}) } \}
\int_{\Delta_w} \omega_w = \\
&((\lambda,\lambda) -(\rho, \rho))\int_{\Delta_w} \omega_w
\endalign
$$

as well as to the whole hypergeometric  system of differential
equations
 and thus provide the basis for
the solutions of this system for generic $\lambda $ and $k$.

\endproclaim

\subhead {Proof of theorem 6.3} \endsubhead
Throughout the proof we do not pay attention to the phases of
multivalued forms, since they do not matter for the differential equations.

The theorem is proved as follows. Consider the following integral
(which is a part of the multiple integral)

$$\int t_{11}^{\lambda_n- \lambda_{n-1} -k} (t_{12} -t_{11})^{k-1}
(t_{22} -t_{11})^{k-1} dt_{11}.$$ Then it satisfies to the following
differential equation

$$\{ {t_{12} t_{22}{\partial^2 \over{\partial t_{12}\partial t_{22}}}} +
(k-1){t_{12} t_{22} \over{ t_{12}- t_{22}}} (\frac {\partial}{\partial
t_{12}}-\frac {\partial}{\partial  t_{22}}) \}
\int t_{11}^{\lambda_n- \lambda_{n-1} -k} (t_{12} -t_{11})^{k-1}
(t_{22} -t_{11})^{k-1} dt_{11} =0
$$

Now let's integrate this differential operator by parts and rewrite it
in terms of variables of the third row ( one should collect the terms
and use homogeneity condition), i.e.
 $ t_{13}, t_{23}, t_{33} $.

One obtains that the integral

$$
\align
\int \int \int & t_{11}^{\lambda_{n}- \lambda_{n-1} -k}
 (t_{12} t_{22})^{\lambda_{n-1}- \lambda_{n-2} -k}
(t_{11} -t_{12})^{k-1}(t_{11} -t_{22})^{k-1}\\
&\times (t_{12} - t_{22})^{2-2 k}
(t_{13} - t_{12})^{k-1}(t_{13} - t_{22})^{k-1}
(t_{23} - t_{12})^{k-1}(t_{23} - t_{22})^{k-1}\\
&\times(t_{33} - t_{12})^{k-1}(t_{33} - t_{22})^{k-1} dt_{11} dt_{12} d
t_{22}
\endalign
$$

satisfies to the following differential equation:

$$
\align
&
\{ {t_{13} t_{23}{\partial^2 \over{\partial t_{13}\partial t_{23}}}}
 +
{t_{13} t_{33}{\partial^2 \over{\partial t_{13}\partial t_{33}}}}
 +
{t_{23} t_{33}{\partial^2 \over{\partial t_{23}\partial t_{33}}}}\\
& +
(k-1)\{{t_{13} t_{23} \over{ t_{13}- t_{23}}} (\frac {\partial}{\partial
t_{13}}-\frac {\partial}{\partial  t_{23}})
+
{t_{23} t_{33} \over{ t_{23}- t_{33}}} (\frac {\partial}{\partial
t_{23}}-\frac {\partial}{\partial  t_{33}})
 +
{t_{13} t_{33} \over{ t_{13}- t_{33}}} (\frac {\partial}{\partial
t_{13}}-\frac {\partial}{\partial  t_{33}}) \} \}\\ &\\
&
\int \int \int  t_{11}^{\lambda_{n}- \lambda_{n-1} -k}
 (t_{12} t_{22})^{\lambda_{n-1}- \lambda_{n-2} -k}
(t_{11} -t_{12})^{k-1}(t_{11} -t_{22})^{k-1}\\
&\times (t_{12} - t_{22})^{2-2 k}
(t_{13} - t_{12})^{k-1}(t_{13} - t_{22})^{k-1}
(t_{23} - t_{12})^{k-1}(t_{23} - t_{22})^{k-1}\\
&\times(t_{33} - t_{12})^{k-1}(t_{33} - t_{22})^{k-1} dt_{11} dt_{12} d
t_{22}
 =\\ &\\
&a_{3}
\int \int \int  t_{11}^{\lambda_{n}- \lambda_{n-1} -k}
 (t_{12} t_{22})^{\lambda_{n-1}- \lambda_{n-2} -k}
(t_{11} -t_{12})^{k-1}(t_{11} -t_{22})^{k-1}\\
&\times (t_{12} - t_{22})^{2-2 k}
(t_{13} - t_{12})^{k-1}(t_{13} - t_{22})^{k-1}
(t_{23} - t_{12})^{k-1}(t_{23} - t_{22})^{k-1}\\
&\times(t_{33} - t_{12})^{k-1}(t_{33} - t_{22})^{k-1} dt_{11} dt_{12} d
t_{22}
\endalign
$$

where $a_{3}=a_{3}(\lambda, k)$ is some constant which takes into
account homogeneity relations. Since finally the eigenvalue will be
determined using Harish-Chandra homomorphism one can avoid calculation of
this constant.

Let's integrate by parts once more and rewrite differential
 operator
in terms of ${t_{i4}, \; i=1,2,3,4}$ and one can do this row
by row.

Calculation (integration by parts) amounts  essentially to the following
identity:

$$
\sum_{i<j, s1,s2}
\frac
{t_{i,m} t_{j,m}}
{(t_{s_1,m+1}-t_{i,m})(t_{s_2,m+1}-t_{j,m})}=
\sum_{s1,s2, i<j}
\frac
{t_{s_1,m+1}t_{s_2,m+1} }
{(t_{s_1,m+1}-t_{i,m})(t_{s_2,m+1}-t_{j,m})}  + C(m)
$$
 where $C(m) =C(m,\lambda,k)$ a constant which takes into account
homogeneity relations.

On the last step the following lemma is used
\proclaim{Lemma 6.4 }\endproclaim
$$
\align
&
\sum_{i<j} ({t_{i,n+1} t_{j,n+1} \frac{\partial^2}{\partial t_{i,n+1} \partial
t_{j,n+1}}}
+(k-1)\frac{ t_{i,n+1} t_{j,n+1} } {t_{i,n+1}-t_{j,n+1}} ({\partial
\over {\partial t_{i,n+1}}}-{\partial \over {\partial t_{j,n+1}}}))
\prod_{p<q}(t_{p,n+1}-t_{q,n+1})^{2 k-1} =\\
&(2 k -1){(n-1)n(n+1)\over 24} (k 6 n -3n +2) \quad  \prod_{p<q} (t_{p,n+1}
-t_{q,n+1})^{2 k-1}
\endalign
$$

Lemma 6.4 easily follows from the next lemma 6.5.

\proclaim{Lemma 6.5}\endproclaim
$$
\align
&
\sum_{i<j} ({t_{i,n+1} t_{j,n+1} \frac{\partial^2}{\partial t_{i,n+1}
\partial t_{j,n+1}}}
)
\prod_{p<q}(t_{p,n+1}-t_{q,n+1}) =\\
&{(n-1)n(n+1)\over 24} (3n +2) \quad  \prod_{p<q} (t_{p,n+1}
-t_{q,n+1})
\endalign
$$
\smallskip
\subhead{Proof of Lemma 6.5}\endsubhead

In fact,

$$
\prod_{p<q}(t_{p,n+1}-t_{q,n+1}) = \sum_{w \in S_{n+1}}
 {(-1)^{n(n+1)\over 2} det(w) \prod t_{w(i),n+1}^{i-1}}
$$

 and lemma follows from the identity:

$$
\sum_{q=2}^{n+1} \sum_{p<q} (p-1)(q-1)={(n-1)n(n+1) (3n +2)\over 24}
$$

\vskip 1cm

\subhead{ Continuation of proof of theorem 6.3 }\endsubhead

$$
\align
\{\sum_{i=1}^{n+1} { ({{t_{i,n+1}{\partial\over{ \partial t_{i,n+1} }
}}})^2}  - k \sum_{i <j} &{
{\frac{t_{j,n+1}+t_{i,n+1}}{t_{j,n+1}-t_{i,n+1}}}
 (t_{i,n+1} {\partial \over
\partial t_{i,n+1}}- t_{j,n+1} {\partial \over \partial t_{j,n+1}}) } \}
\int_{\Delta_w} \omega_w = \\
& c_n \int_{\Delta_w} \omega_w
\endalign
$$

Finally, the eigenvalue can be determined using Harish-Chandra
homomorphism

$$
c_n=(w\lambda +\rho, w\lambda +\rho) - 2(\rho, w \lambda +\rho)=
(\lambda, \lambda) -(\rho, \rho)
$$
So

$$
\align
\{\sum_{i=1}^{n+1} { ({{t_{i,n+1}{\partial\over{ \partial t_{i,n+1} }
}}})^2}  - k \sum_{i <j} &{
{\frac{t_{j,n+1}+t_{i,n+1}}{t_{j,n+1}-t_{i,n+1}}}
 (t_{i,n+1} {\partial \over
\partial t_{i,n+1}}- t_{j,n+1} {\partial \over \partial t_{j,n+1}}) } \}
\int_{\Delta_w} \omega_w = \\
&((\lambda,\lambda) -(\rho, \rho))\int_{\Delta_w} \omega_w
\endalign
$$

This  completes the proof that the integrals satisfy to the second
order differential equation and thus to the whole hypergeometric
system of differential equations .

\smallskip

\remark{Remark 6.6} This elementary proof does not use the advanced theory of
Knizhnik-Zamolodchikov equation (as well as Dunkl operators)
,but uses only integration by parts and the ability to present
asymptotic solutions.
Though one should notice that integration by parts is a commonly used
technique in working with similar integrals appearing in the theory of
Knizhnik-Zamolodchikov equation.
\endremark
\vskip 1cm

Let $F_w$ be the normalized asymptotic solution of Heckman-Opdam
hypergeometric system of type $A_n$, i.e.
 $F_w(z) = z^{w \lambda +\rho}(1 +\ldots)$.

Now let $ z_{1}, z_{2} ,\ldots,  z_{n+1}$
approach to $1$ , while keeping inequalities

$$0 < | z_{1}| < |z_{2}| <\ldots < |z_{n+1}|. $$

\proclaim{Theorem 6.7}(Opdam)

$$
F_w(1)=
\lim_{z \rightarrow 1} {F_{w}(z)}=
{\frac
{\prod_{\alpha \in R_{+}}
      {\frac {\Gamma( (w \lambda, \alpha^{\vee})+1)}
      {\Gamma( (w \lambda, \alpha^{\vee})- k +1)}}}
{ \prod_{\alpha \in R_{+}}
 {\frac
 {\Gamma( -(\rho , \alpha^{\vee})+1)}
      {\Gamma( -( \rho, \alpha^{\vee})- k +1)}}}}
$$

cf. theorem 6.3 [21].
\endproclaim

\proclaim {Theorem 6.8} The limit of integral of $\omega_w$ over
$\Delta_w$ as all $z_{i}$ approach to $1$ , while preserving the
above inequalities, is equal to :

$$
\align
 & \lim_{z \rightarrow 1} \int_{\Delta_w(z)} \omega_w= \\
&\prod_{\alpha \in{R_{+}}} \sin(\pi ((-w\lambda,\alpha^\vee ) +k)) \times
e^{-2 \pi i(\lambda, \delta)} e^{-\pi i (k-1) l(w)} (2 i)^{\frac{n(n+1)}2}\\
&{\sin(\pi k)^{n+1}\over{{\sin (\pi k)} { \sin(2 \pi k) }\ldots {\sin((n+1)
\pi k)}}}\times {\Gamma(k)^{(n+1)(n+2)\over2}\over {\Gamma(k) \Gamma(2 k)
\ldots
\Gamma((n+1)k)}}
\endalign
$$
\endproclaim
The theorem easily follows from theorems 6.1 , 6.7, 6.6.
see also [21, 17, 33].

\subhead Acknowledgments \endsubhead

I would like to thank I.~Gelfand for stimulating discussions
concerning   the theory of hypergeometric
functions and the theory of spherical functions and  for the
suggestions on the organization of the material of the paper.
I would like to thank S.~Lukyanov for stimulating discussions concerning
bosonization technique in conformal field theory.
 I would also
like to thank V.~Brazhnikov, V.~Dolotin, D.~Fradkin for helpful discussions and
A.~Guglin for the help in preparation of the manuscript.

\midinsert\vskip 7cm
\includegraphics{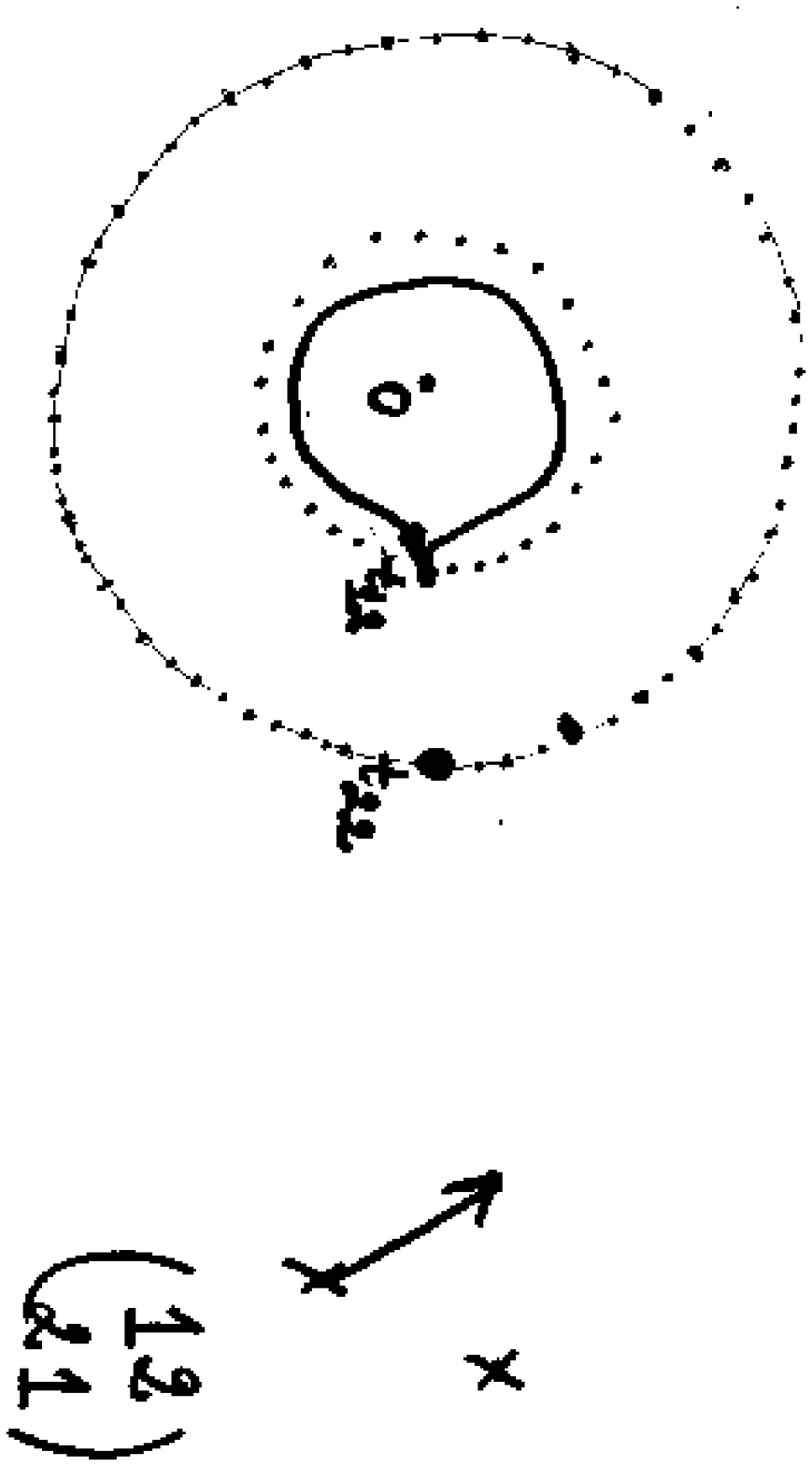}
\botcaption{Figure 4a}
\endcaption
\endinsert

\midinsert\vskip 7cm
\includegraphics{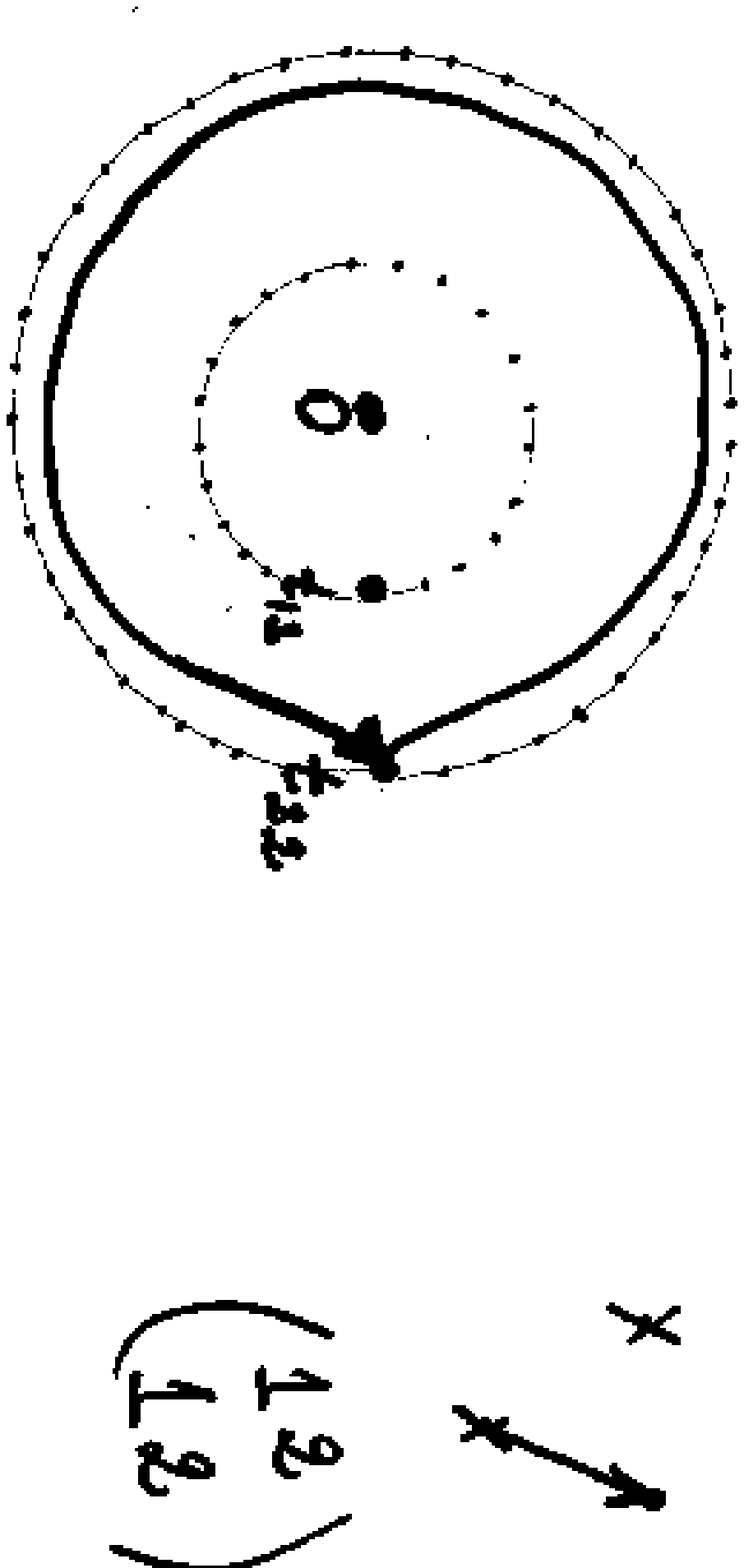}
\botcaption{Figure 4b}
\endcaption
\endinsert

\midinsert\vskip 7cm
\includegraphics{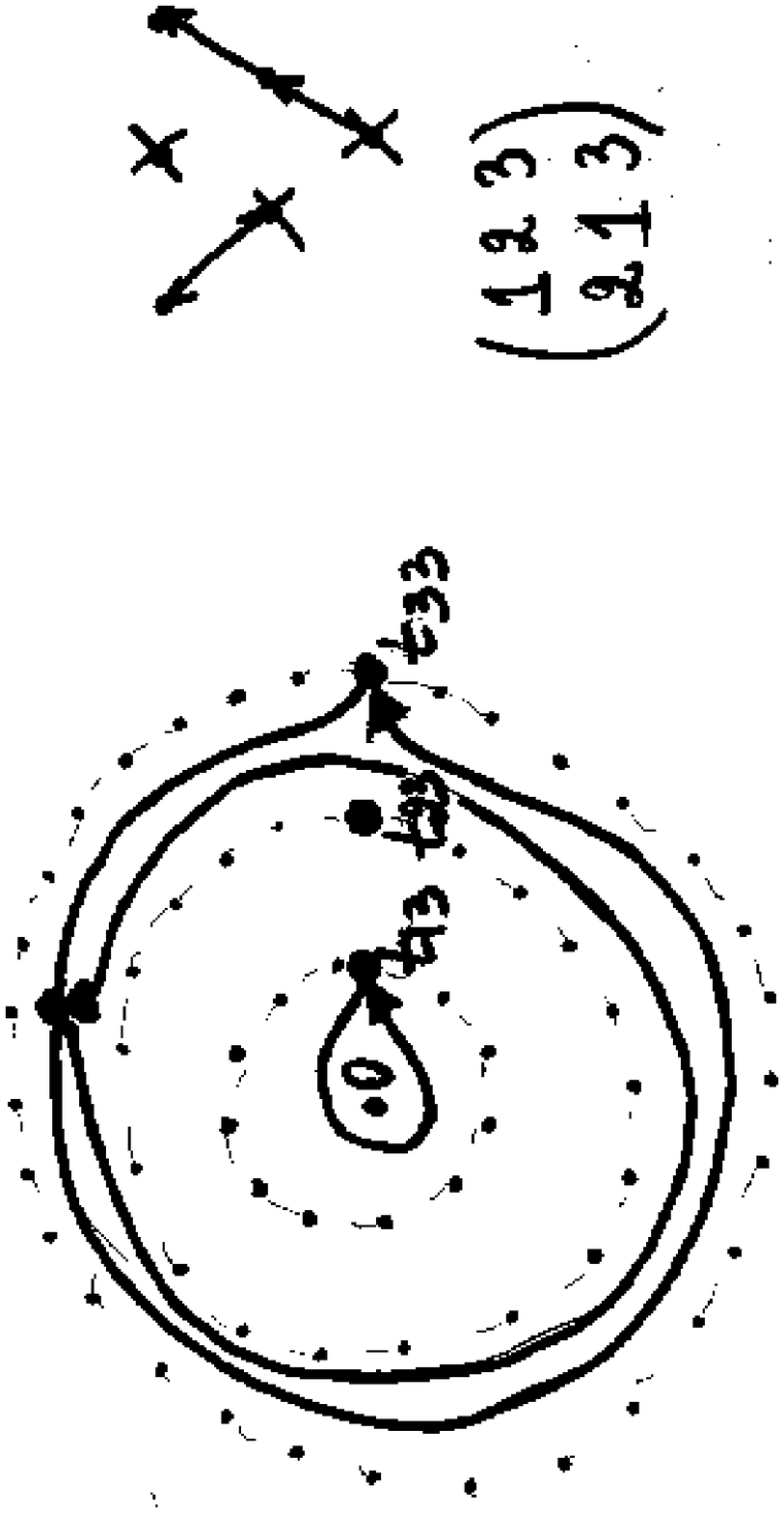}
\botcaption{Figure 5a}
\endcaption
\endinsert

\midinsert\vskip 7cm
\includegraphics{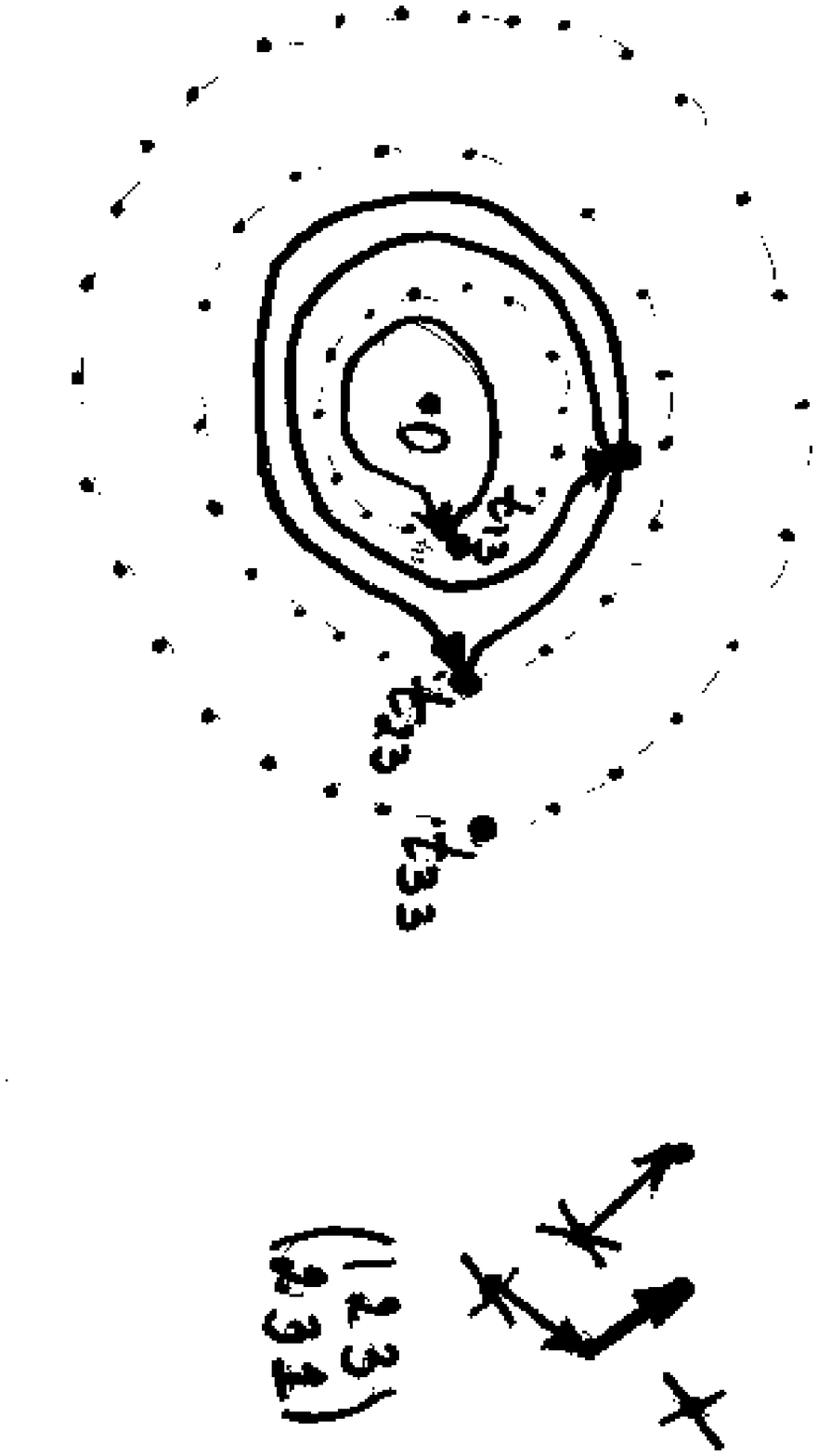}
\botcaption{Figure 5b}
\endcaption
\endinsert

\midinsert\vskip 7cm
\includegraphics{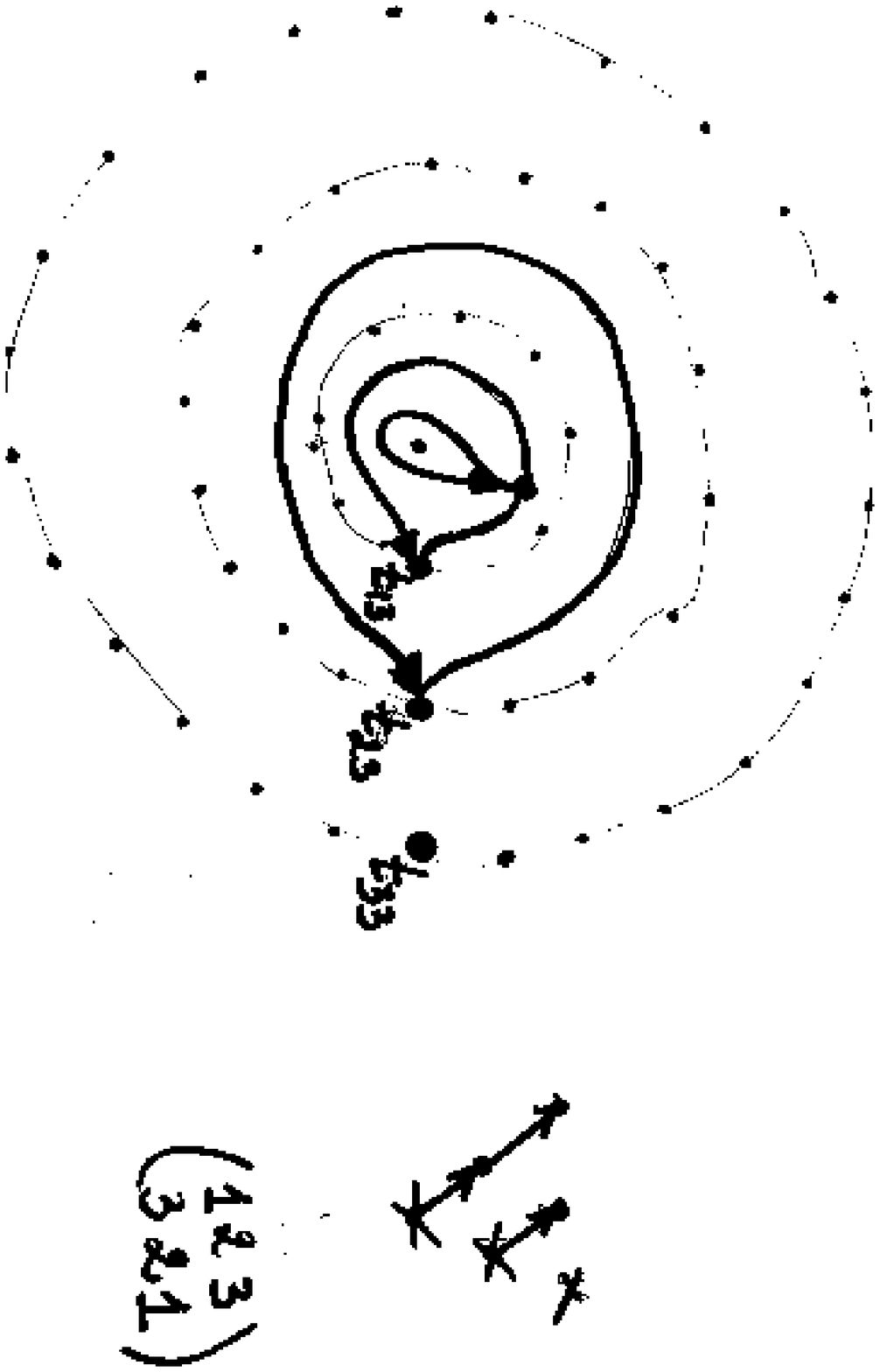}
\botcaption{Figure 5c}
\endcaption
\endinsert

\midinsert\vskip 7cm
\includegraphics{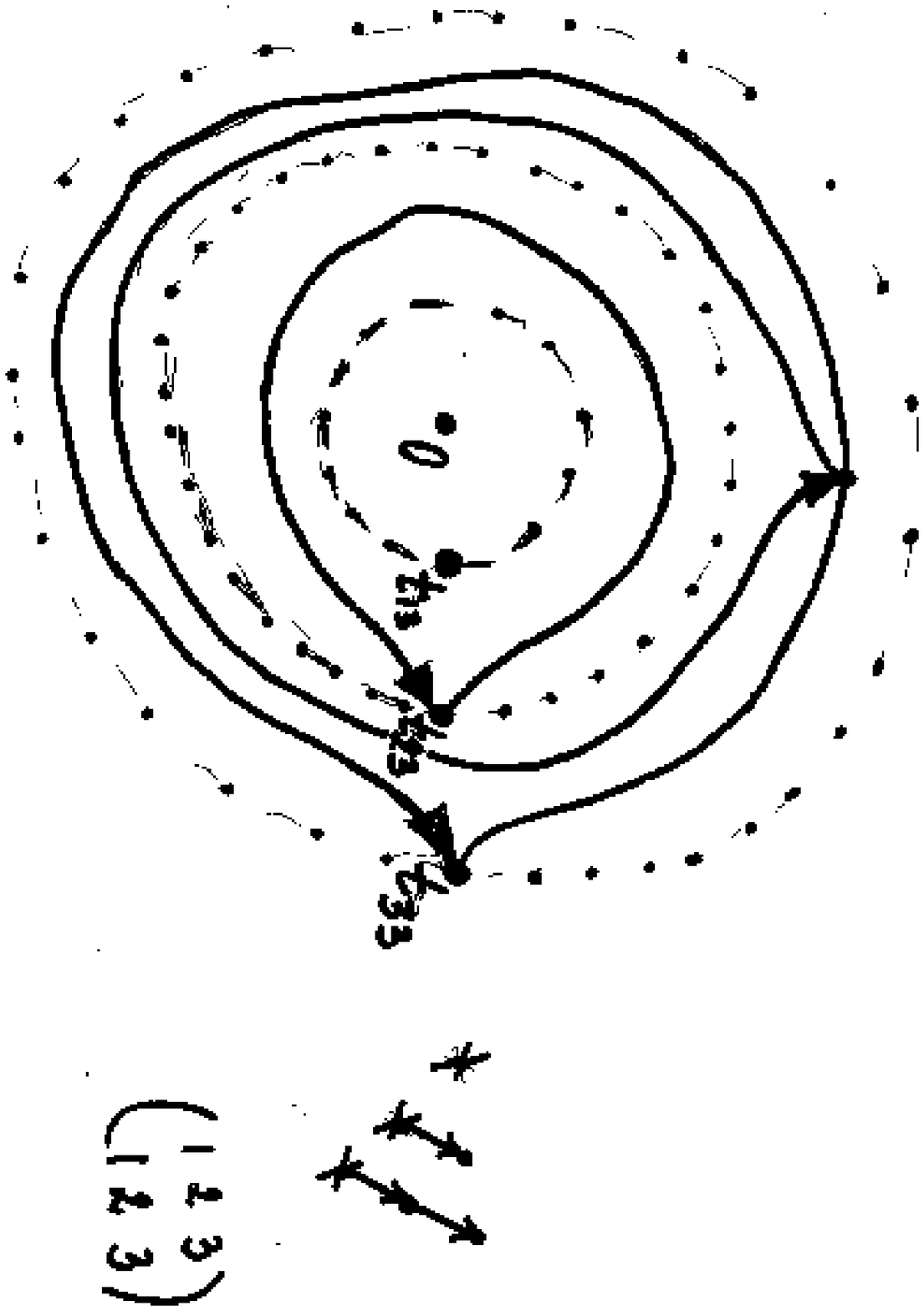}
\botcaption{Figure 5d}
\endcaption
\endinsert

\midinsert\vskip 7cm
\includegraphics{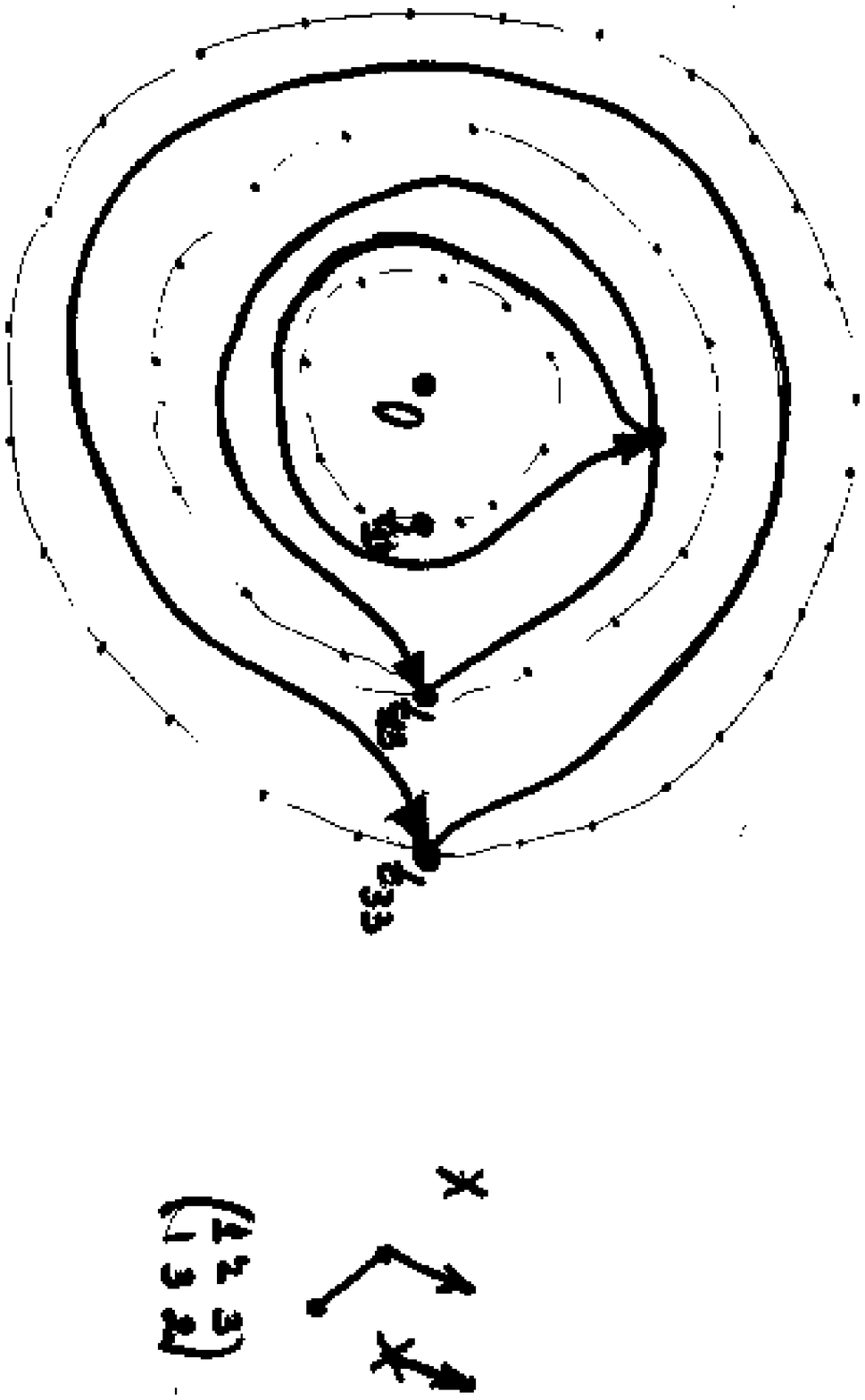}
\botcaption{Figure 5e}
\endcaption
\endinsert

\midinsert\vskip 7cm
\includegraphics{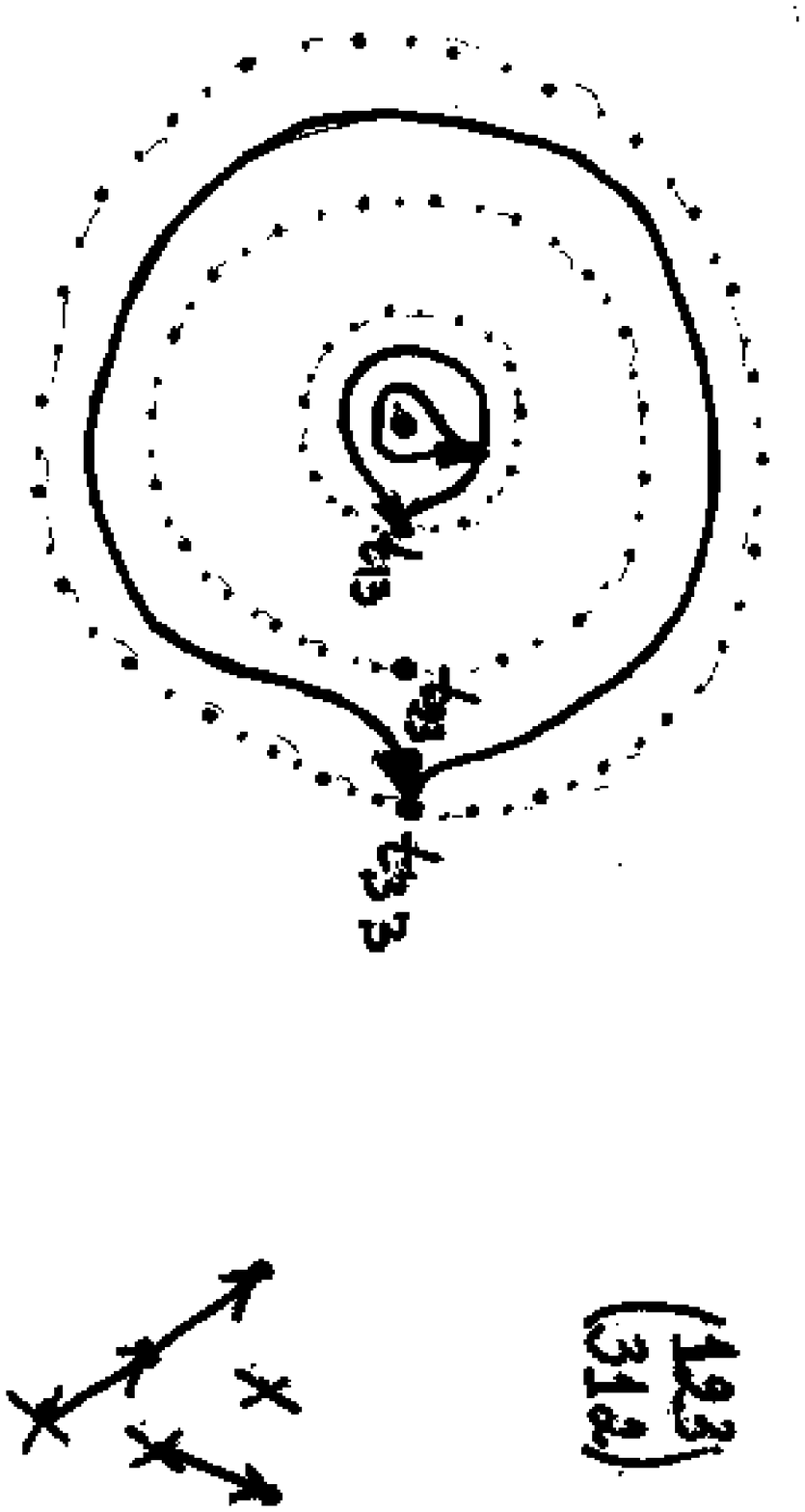}
\botcaption{Figure 5f}
\endcaption
\endinsert

\enddocument